\documentclass{amsart}
\usepackage{amssymb}

\newtheorem{theorem}{Theorem}[section]
\newtheorem{lemma}[theorem]{Lemma}
\numberwithin{equation}{section}
\begin{document}

\title{$(k,r)$-admissible configurations and intertwining operators}
\author{Mirko Primc}
\address{Department of Mathematics, University of Zagreb, Bijeni\v{c}ka 30, Zagreb, Croatia}
\email{primc@math.hr}
\thanks{Partially supported  by the Ministry of Science and Technology of the
Republic of Croatia, grant 0037125.} 

\begin{abstract}
Certain combinatorial bases of Feigin-Stoyanovsky's type subspaces
of level $k$ standard modules for affine Lie algebra
$\mathfrak{sl}(r,\mathbb C)\,\widetilde{}\ $ are parametrized by
$(k,r)$-admissible configurations. In this note we use
Capparelli-Lepowsky-Milas' method to give a new proof of linear
independence of these bases, the main ingredient in the proof being
the use of Dong-Lepowsky's intertwining operators for fundamental
$\mathfrak{sl}(r,\mathbb C)\,\widetilde{}\ $-modules.
\end{abstract}
\maketitle

\section{Introduction}
Let $\mathfrak g$ be a simple Lie algebra, $\mathfrak h$ a Cartan
subalgebra of $\mathfrak g$ and $\mathfrak g=\mathfrak
h+\sum_{\alpha\in R}\mathfrak g_\alpha$ the root space decomposition
with fixed root vectors $x_\alpha$. Let
\begin{equation}\label{decomposition of g}
\mathfrak g={\mathfrak g}_{-1}+{\mathfrak g}_{0}+{\mathfrak g}_{1}
\end{equation}
be a $\mathbb Z$-grading of $\mathfrak g$ such that ${\mathfrak
h}\subset{\mathfrak g}_{0}$. Let $\tilde{\mathfrak g} = {\mathfrak
g}\otimes \mathbb C [t,t\sp{-1}]+\mathbb C c+\mathbb C d$ be the
associated affine Lie algebra with the canonical central element $c$
and fixed real root vectors $x_\alpha(n)=x_\alpha\otimes t\sp n$.
Let
\begin{equation*}
\tilde{\mathfrak g}=\tilde{\mathfrak g}_{-1}+\tilde{\mathfrak
g}_{0}+\tilde{\mathfrak g}_{1}
\end{equation*}
be the decomposition corresponding to (\ref{decomposition of g}),
where $\tilde{\mathfrak g}_{1}={\mathfrak g}_{1}\otimes \mathbb C
[t,t\sp{-1}]$ is a commutative Lie subalgebra with a basis
\begin{equation}\label{a basis of g1}
\{x_\gamma(j)\mid j\in\mathbb Z, \gamma\in\Gamma\}
\end{equation}
for some $\Gamma\subset R$. Let $L(\Lambda)$ be a standard (i.e.,
integrable highest weight) $\tilde{\mathfrak g}$-module of level
$k=\Lambda(c)$ with a fixed highest weight vector $v_\Lambda$.
Define a Feigin-Stoyanovsky's type subspace $W(\Lambda)$ as
$$
W(\Lambda)=U(\tilde{\mathfrak g}_{1})v_\Lambda.
$$
This subspace of $L(\Lambda)$ is analogous to a principal subspace
introduced in \cite{FS}.

It seems that a spanning set of $W(\Lambda)$ consisting of
``monomial vectors'' $x(\pi)v_\Lambda$, where $x(\pi)$ are monomials
in basis elements (\ref{a basis of g1}), can be reduced to a basis
of $W(\Lambda)$ by using ${\mathfrak g}_{0}$-modules of relations on
$W(\Lambda)\subset L(\Lambda)$ generated by the adjoint action of
${\mathfrak g}_{0}$ on relations
\begin{equation}\label{relations}
\sum_{j_1+\dots+j_{k+1}=n}x_\gamma(j_1)\dots
x_\gamma(j_{k+1})=0,\quad n\in\mathbb Z,
\end{equation}
for some long root $\gamma\in\Gamma$. For $\mathfrak
g=\mathfrak{sl}(r,\mathbb C)$ and a particular choice of
(\ref{decomposition of g}) such bases are constructed in \cite{P1}
for all $W(\Lambda)$, and it turns out that basis elements are
parametrized by $(k,r)$-admissible configurations --- combinatorial
objects introduced and studied in a series of papers
\cite{FJLMM}--\cite{FJMMT}. On the other side, for any classical
simple Lie algebra $\mathfrak g$ and any choice of
(\ref{decomposition of g}) such bases are constructed in \cite{P2},
but only for $W(\Lambda_0)$
--- a subspace of the basic $\tilde{\mathfrak g}$-module
$L(\Lambda_0)$.

The problem of constructing monomial bases of $W(\Lambda)$ is a part
of Lepowsky-Wilson's approach to Rogers-Ramanujan type identities
initiated in \cite{LW}, and the choice of vertex operator algebra
relations (\ref{relations}) is close to constructions in \cite{LP}
and \cite{MP}. In Lepowsky-Wilson's approach a lack of known
combinatorial identities usually makes a proof of linear
independence of combinatorial basis a hard problem. In \cite{P2}
linear independence of monomial basis of $W(\Lambda_0)$ is proved by
using the crystal base character formula \cite{KKMMNN}, but it is
not clear how to extend the proof for levels $k>1$. On the other
side, in \cite{P1} linear independence is proved directly by using
``Schur functions'', but it is not clear how to extend the proof to
other simple Lie algebras $\mathfrak g$.

G.~Georgiev proves in \cite{G} linear independence of
quasi-particle bases of Feigin-Stoyanovsky's level $k$ principal
subspaces for $\mathfrak{sl}(r,\mathbb C)\,\widetilde{}\ $ by
using Dong-Lepowsky's intertwining operators. Also by using
intertwining operators, in \cite{CLM1} and \cite{CLM2}
Rogers-Selberg recursions for characters of level $k$  principal
subspaces for $\mathfrak{sl}(2,\mathbb C)\,\widetilde{}\ $ are
obtained. Although S. Capparelli, J. Lepowsky and A. Milas  avoid
the explicit use of combinatorial bases, their way of using
intertwining operators leads to a simple proof of linear
independence of the underlying monomial bases for all $k$. Here we
present the level $1$ case when the proof coincides with
Georgiev's argument: For $\mathfrak g=\mathfrak{sl}(2,\mathbb
 C)=\text{span}\,\{x_\alpha, h, x_{-\alpha}\}$ bases of principal
subspaces $W(\Lambda_0)$ and $W(\Lambda_1)$ of fundamental
$\tilde{\mathfrak g}$-modules consist of monomial vectors (for
$i=0,1$)
\begin{equation}\label{sl2 bases}
x(\pi)v_{\Lambda_i}=\prod_{j\geq1}x_\alpha(-j)\sp{a_{j-1}}v_{\Lambda_i}
\end{equation}
such that $\pi=(a_0,a_1,\dots)$ are $(1,2)$-admissible
configurations, i.e., such that
$$
a_j+a_{j-1}\leq 1\quad\text{for all}\quad j\geq 1,\quad\text{}\quad
a_0\leq 1-i,
$$
with finitely many non-zero $a_j$. The key role in the proof of
linear independence of (\ref{sl2 bases}) play two maps constructed
from the intertwining operators, $[1]$ and $[\omega]$ in our
notation below,
\begin{align*}
&[1]\colon W(\Lambda_0)\to W(\Lambda_1),\quad
[1]v_{\Lambda_0}=v_{\Lambda_1},\\
&[\omega]\colon W(\Lambda_0)\to W(\Lambda_1),\quad
[\omega]v_{\Lambda_0}=v_{\Lambda_1},\\
&[\omega]\colon W(\Lambda_1)\to W(\Lambda_0),\quad
[\omega]v_{\Lambda_1}=x_\alpha(-1)v_{\Lambda_0},
\end{align*}
such that $[\omega]$'s are injective, and such that for all
$j\in\mathbb Z$
$$
x_\alpha(j)[1]=[1]x_\alpha(j),\quad
x_\alpha(j)[\omega]=[\omega]x_\alpha(j+1).
$$
The proof goes by induction on degree $n=\sum ja_{j-1}$ of monomials
$x(\pi)$. First consider a linear combination $\sum c_\pi
x(\pi)v_{\Lambda_1}=0$ of vectors (\ref{sl2 bases}) of degree $n$ in
$W(\Lambda_1)$. Then
$$
0=\sum c_\pi\left(\prod_{j\geq
2}x_\alpha(-j)\sp{a_{j-1}}\right)[\omega]v_{\Lambda_0}=[\omega]\sum
c_\pi\left(\prod_{j\geq
2}x_\alpha(-j+1)\sp{a_{j-1}}\right)v_{\Lambda_0}.
$$
Since $[\omega]$ is injective, and vectors $x(\pi\sp+)v_{\Lambda_0}$
on the right hand side are $(1,2)$-admissible of degree $<n$, by
induction hypothesis we get $c_\pi=0$. Now consider a linear
combination $\sum c_\pi x(\pi)v_{\Lambda_0}=0$ of vectors (\ref{sl2
bases}) of degree $n$ in $W(\Lambda_0)$. Then
$$
0=\sum c_\pi x(\pi)v_{\Lambda_0}\overset{[1]}\longrightarrow
\sum{}\sp{{}'} c_\pi x(\pi)v_{\Lambda_1}=0,
$$
where $\sum'$ runs over all $\pi$ such that $a_0=0$, and, by already
proven independence in $W(\Lambda_1)$, we get $c_\pi=0$ for all such
$\pi$. So we have $0=\sum c_\pi x(\pi)v_{\Lambda_0}=\sum '' c_\pi
x(\pi)v_{\Lambda_0}$, where $\sum ''$ runs over all $\pi$ such that
$a_0=1$, that is
$$
0=\sum c_\pi\left(\prod_{j\geq
2}x_\alpha(-j)\sp{a_{j-1}}\right)x_\alpha(-1)v_{\Lambda_0}=\sum
c_\pi\left(\prod_{j\geq
2}x_\alpha(-j)\sp{a_{j-1}}\right)[\omega]v_{\Lambda_1}.
$$
Now we commute $[\omega]$ to the left and, as before, we conclude
that $c_\pi=0$ for all $\pi$.

In this note we extend this proof for $\mathfrak{sl}(2,\mathbb
C)\,\widetilde{}\ $ level $1$ to a proof of linear independence of
all combinatorial bases constructed in \cite{P1}. One should hope
that this new proof of linear independence of combinatorial bases
can be extended further to all Feigin-Stoyanovsky's type subspaces
for all classical affine Lie algebras.

I thank Haisheng Li for many useful discussions and for pointing out
the connection between simple currents and the ``constant factor''
$[\omega]=e\sp{\alpha/2}$ in \cite{CLM1}.

\section{Affine Lie algebra $\mathfrak{sl}(\ell +1,\mathbb C)\,\widetilde{}\ $}

From now on let ${\mathfrak g}={\mathfrak sl}(\ell+1,\mathbb C)$ and
let $\mathfrak h$ be the Cartan subalgebra of diagonal matrices.
Denote, as usual, the corresponding root system
$$
R=\{\pm(\varepsilon_i-\varepsilon_j)\mid 1\leq i<j\leq \ell+1\},
$$
fix simple roots $\alpha_1=\varepsilon_1-\varepsilon_2$, \dots,
$\alpha_\ell=\varepsilon_\ell-\varepsilon_{\ell+1}$ and denote by
$\omega_1$, \dots, $\omega_\ell$ the corresponding fundamental
weights of ${\mathfrak g}$. It will be convenient to write
$\omega_0=0$. Set
$$
\Gamma=\{\gamma_1,\gamma_2,\dots,\gamma_\ell\},\quad
\gamma_i=\varepsilon_i-\varepsilon_{\ell+1}=\alpha_i+\dots+\alpha_\ell.
$$
Note that $\gamma_1=\theta$ is the maximal root and that $\Gamma$ is
a basis of  $\mathfrak h\sp*$. Denote by $Q=Q(R)$ the root lattice
and by $P=P(R)$ the weight lattice of $\mathfrak g$. Denote by
$\langle\cdot,\cdot\rangle$ the normalized Killing form such that
$\langle\theta,\theta\rangle=2$, where we identify $\mathfrak h\cong
\mathfrak h\sp*$ via $\langle\cdot,\cdot\rangle$. For each root
$\alpha$ fix a root vector $x_\alpha$.

Denote by $\tilde{\mathfrak g}$ the associated affine Lie algebra
$$
\tilde{\mathfrak g} =\coprod_{n\in\mathbb Z}{\mathfrak g}\otimes
t^{n}+\mathbb C c+\mathbb C d
$$
with the canonical central element $c$ and the degree element $d$
such that $[d,x\otimes t^{n}]=-n\,x\otimes t^{n}$ (cf. \cite{K}).
Denote by $\Lambda_0$, $\Lambda_1$, \dots, $\Lambda_\ell$ the
corresponding fundamental weights of $\tilde{\mathfrak g}$. Write
$x(n)=x\otimes t^{n}$ \ for $x\in{\mathfrak g}$ and $n\in\mathbb Z$
and denote by $x(z)=\sum_{n\in\mathbb Z} x(n) z^{-n-1}$  a formal
Laurent series in formal variable $z$.

\section{Feigin-Stoyanovsky's type subspaces $W(\Lambda)$}

Denote by $L(\Lambda)$ a standard $\tilde{\mathfrak g}$-module with
a dominant integral highest weight
$$
\Lambda=k_0\Lambda_0+k_1\Lambda_1+\dots+k_\ell\Lambda_\ell,
$$
$k_0, k_1, \dots, k_\ell\in\mathbb Z_{+}$ (cf. \cite{K}). Denote
by $k=\Lambda(c)$ the level of $\tilde{\mathfrak g}$-module
$L(\Lambda)$,
$$
k=k_0 +k_1+ \dots+ k_\ell.
$$
For each fundamental $\tilde{\mathfrak g}$-module $L(\Lambda_i)$
fix a highest weight vector $v_{\Lambda_i}$. By complete
reducibility of tensor products of standard modules, for level
$k>1$ we have
$$
L(\Lambda)\subset L(\Lambda_{\ell})\sp{\otimes k_\ell}\otimes\dots
\otimes L(\Lambda_{1})\sp{\otimes k_1}\otimes
L(\Lambda_{0})\sp{\otimes k_0}
$$
with the highest weight vector
$$
v_{\Lambda} = v_{\Lambda_{\ell}}\sp{\otimes k_\ell}\otimes\dots
\otimes v_{\Lambda_{1}}\sp{\otimes k_1}\otimes
v_{\Lambda_{0}}\sp{\otimes k_0}.
$$
Set
$$
{\tilde{{\mathfrak g}}_1}=\text{span}\{x_\gamma(n)\mid \gamma \in
\Gamma, n\in\mathbb Z \}.
$$
Note that ${\tilde{{\mathfrak g}}_1}$ is a commutative Lie
subalgebra of $\tilde{{\mathfrak g}}$. For each integral dominant
$\Lambda$ define a Feigin-Stoyanovsky's type subspace
$$
W(\Lambda)=U({\tilde{{\mathfrak g}}_1})v_\Lambda\subset L(\Lambda).
$$

Denote by $\pi\colon \{x_\gamma(-j)\mid \gamma \in \Gamma, j\geq
1\}\to \mathbb Z_+$ \ a ``colored partition'' for which a finite
number of ``parts'' $x_\gamma(-j)$ (of degree $j$ and color
$\gamma$) appear $\pi(x_\gamma(-j))$ times, and denote by
$$
x(\pi)=\prod x_\gamma(-j)\sp{\pi(x_\gamma(-j))}\in
U({\tilde{{\mathfrak g}}_1})=S({\tilde{{\mathfrak g}}_1})
$$
the corresponding monomials. We can identify $\pi$ with a sequence
$(a_i)_{i=0}\sp\infty$ with finitely many non-zero terms
$a_{\ell(j-1)+r-1}=\pi(x_{\gamma_r}(-j))$ and
$$
x(\pi)=\dots\,
x_{\gamma_1}(-2)\sp{a_\ell}x_{\gamma_\ell}(-1)\sp{a_{\ell-1}}\dots
x_{\gamma_1}(-1)\sp{a_0}.
$$
By Poincar\' e-Birkhoff-Witt theorem we have a spanning set of
monomial vectors $x(\pi)v_\Lambda$ in $W(\Lambda)$. By using VOA
relations
$$
x_{\beta_{1}}(z)\cdots x_{\beta_{k+1}}(z)=0,
\quad\quad\beta_{1},\dots, \beta_{k+1}\in\Gamma,
$$
for level $k$ standard $\tilde{\mathfrak g}$-modules, we can reduce
PBW spanning set to a spanning set of $(k,\ell+1)$-admissible
monomial vectors
$$
\dots
x_{\gamma_1}(-2)\sp{a_\ell}x_{\gamma_\ell}(-1)\sp{a_{\ell-1}}\dots
x_{\gamma_1}(-1)\sp{a_0}v_\Lambda
$$
satisfying difference conditions
\begin{equation}\label{admissible difference conditions}
0\leq a_i\leq k,\quad a_i+\dots+a_{i+\ell}\leq k
\end{equation}
for all $i\in\mathbb Z_+$, and initial conditions
\begin{equation}\label{admissible initial conditions}
a_0\leq k_0, \quad a_0+a_{1}\leq k_0+k_1, \quad \dots \quad
a_0+\dots+a_{\ell-1}\leq k_0+\dots+k_{\ell-1}.
\end{equation}
This spanning set is a basis (cf. \cite{FJLMM} and \cite{P1}). {\it
As said before, the aim of this note is to give a new proof of
linear independence of this spanning set, and hence a new proof of
the following:}
\begin{theorem}
The set of $(k,\ell+1)$-admissible monomial vectors
$x(\pi)v_\Lambda$ is a basis of $W(\Lambda)$.
\end{theorem}

It should be noted that the present formulation of initial and
difference conditions for basis elements $x(\pi)v_\Lambda$, due to
\cite{FJLMM}, is crucial for arguments which follow --- see
(\ref{diff.cond.-init.cond.1}) and (\ref{diff.cond.-init.cond.2})
below. By following \cite{FJLMM}, define $(k,\ell+1)$-admissible
configurations as sequences $(a_i)_{i=0}\sp\infty$ with finitely
many non-zero terms satisfying initial and difference conditions
(\ref{admissible difference conditions})--(\ref{admissible initial
conditions}).

\section{Vertex operator construction of level $1$ modules}

We use the well known Frenkel-Kac-Segal construction (\cite{FK},
\cite{S}, cf. \cite{FLM}) of fundamental $\tilde{\mathfrak
g}$-modules $L(\Lambda_i)$ on the tensor product $M(1)\otimes
{\mathbb C}[P]$ of the Fock space $M(1)$ for the homogeneous
Heisenberg subalgebra and the group algebra ${\mathbb C}[P]$ of
the weight lattice with a basis $e\sp\lambda$, $\lambda\in P$. The
action of Heisenberg subalgebra on $M(1)\otimes {\mathbb C}[P]$
extends to the action of Lie algebra $\tilde{\mathfrak g}$ via the
vertex operator formula
\begin{equation}\label{vertex operator formula}
x_\alpha(z)=E\sp-(-\alpha,z)E\sp+(-\alpha,z) e\sb\alpha z\sp\alpha
\end{equation}
for properly chosen root vectors $x_\alpha$, where
$z\sp\alpha=1\otimes z\sp\alpha$, $z\sp\alpha
e\sp\lambda=z\sp{\langle\alpha,\lambda\rangle}$, and
$$
E^\pm(\alpha,z) = E^\pm(\alpha,z)\otimes 1 =\exp \biggl(\sum_{n>0}
\alpha (\pm n) z^{\mp n}\big/ (\pm n)\biggr)\otimes 1.
$$
Then
$$
M(1)\otimes {\mathbb
C}[P]=L(\Lambda_0)+L(\Lambda_1)+\dots+L(\Lambda_\ell)
$$
as $\tilde{\mathfrak g}$-module, and for $i=0,1,\dots,\ell$ we fix
the highest weight vectors
$$
v_{\Lambda_i}=1\otimes e\sp{\omega_i}\quad\text{in}\quad
L(\Lambda_i)=M(1)\otimes e\sp{\omega_i}{\mathbb C}[Q].
$$

Since we also use Dong-Lepowsky's level $1$ intertwining operators
$\mathcal Y$, or, to be more precise, since we use some coefficients
of operators
\begin{equation}\label{intertwining operators}
{\mathcal Y}(1\otimes e\sp\lambda,z)=
E\sp-(-\lambda,z)E\sp+(-\lambda,z) e\sb\lambda z\sp\lambda
e\sp{i\pi\lambda}c(\cdot,\lambda)
\end{equation}
for $\lambda\in P$, we adopt the construction and notation from
[DL, equation (12.3)], where $e_\lambda=1\otimes e_\lambda=
1\otimes e\sp\lambda\epsilon(\lambda,\cdot)$ is defined in [DL,
equations (13.1) and (13.6)] and $c(\alpha,\beta)$ is defined in
[DL, equation (12.52)].

\section{Initial conditions for standard modules}

Since
$\langle\gamma_i,\omega_j\rangle=\langle\alpha_i+\dots+\alpha_\ell,\omega_j\rangle$,
it follows from the vertex operator formula (\ref{vertex operator
formula}) that
\begin{equation}\label{level 1 initial conditions}
x_{\gamma_i}(-1)v_{\Lambda_j}\neq 0\quad \text{if and only if
}\quad i >j.
\end{equation}
In particular
\begin{align*}
&x_{\gamma_1}(-1)v_{\Lambda_1}=x_{\gamma_1}(-1)v_{\Lambda_2}=\dots
x_{\gamma_1}(-1)v_{\Lambda_\ell}=0,\\
&x_{\gamma_2}(-1)v_{\Lambda_2}=x_{\gamma_2}(-1)v_{\Lambda_3}=\dots
x_{\gamma_2}(-1)v_{\Lambda_\ell}=0,
\end{align*}
etc, so we could say that $x_{\gamma_1}(-1)$ can ``nontrivially
sit'' only on $v_{\Lambda_0}$, that $x_{\gamma_2}(-1)$ can
``nontrivially sit'' only on $v_{\Lambda_0}$ and $v_{\Lambda_1}$,
and so on. Moreover, for $\alpha,\beta\in\Gamma$ we have a VOA
relation $x_\alpha(z)x_\beta(z)=0$ on every level $1$ standard
module (cf. \cite{P1}), which  implies  ``Pauli's exclusion
principle"
$$
x_\alpha(-1)x_\beta(-1)v_{\Lambda_j}=0\quad\text{for}\quad
\alpha,\beta\in\Gamma, \ j=0,\dots,\ell.
$$

Now consider level $k>1$ standard modules; we view $L(\Lambda)$
embedded in the tensor product of fundamental modules with the
highest weight vector
$$
v_{\Lambda} = v_{\Lambda_{\ell}}\sp{\otimes k_\ell}\otimes\dots
\otimes v_{\Lambda_{1}}\sp{\otimes k_1}\otimes
v_{\Lambda_{0}}\sp{\otimes k_0}.
$$
Then by (\ref{level 1 initial conditions}), Pauli's exclusion
principle and Dirichlet's pigeonhole principle we see that $a_0$
factors of $x_{\gamma_1}(-1)\sp{a_0}$ can nontrivially sit on $k_0$
factors of $v_{\Lambda_{0}}\sp{\otimes k_0}$ only if $a_0\leq k_0$,
and $a_0+a_1$ factors of
$x_{\gamma_2}(-1)\sp{a_1}x_{\gamma_1}(-1)\sp{a_0}$ can nontrivially
sit on $k_0+k_1$ factors of $v_{\Lambda_{1}}\sp{\otimes k_1}\otimes
v_{\Lambda_{0}}\sp{\otimes k_0}$ only if $a_0+a_{1}\leq k_0+k_1$,
etc. In this way we see that
\begin{equation}\label{level k initial conditions for vectors}
x_{\gamma_\ell}(-1)\sp{a_{\ell-1}}\dots
x_{\gamma_1}(-1)\sp{a_0}v_\Lambda\neq 0
\end{equation}
if and only if
\begin{equation}\label{level k initial conditions}
a_0\leq k_0, \quad a_0+a_{1}\leq k_0+k_1, \quad \dots \quad
a_0+\dots+a_{\ell-1}\leq k_0+\dots+k_{\ell-1}.
\end{equation}
Hence the initial conditions (\ref{admissible initial conditions})
for $(k,\ell+1)$-admissible monomial vectors are equivalent to
(\ref{level k initial conditions for vectors}). Note that monomial
vectors (\ref{level k initial conditions for vectors}) are linearly
independent because different vectors have different $\mathfrak
h$-weights.

For a fixed level $k$ it is convenient to abbreviate the notation:
write
$$
{\mathcal A}=(a_{\ell-1},\dots,a_0)\leq
(k_{\ell-1},\dots,k_0)={\mathcal K}
$$
if (\ref{level k initial conditions}) holds. This is obviously a
partial order on the set of $\ell$-tuples of integers. Write
$x\sp{\mathcal A}=x_{\gamma_\ell}(-1)\sp{a_{\ell-1}}\dots
x_{\gamma_1}(-1)\sp{a_0}$ and $v_{\mathcal K}=v_\Lambda$. Then the
equivalence of (\ref{level k initial conditions for vectors}) and
(\ref{level k initial conditions}) can be restated briefly as
\begin{lemma}\label{L:initial conditions} \quad $x\sp{\mathcal
A}v_{\mathcal K}\ \neq 0$\quad if and only if \quad ${\mathcal
A}\leq{\mathcal K}$.
\end{lemma}

\section{Simple current operator}

As above, for $\lambda\in P$ we consider $e\sp\lambda$ as an
operator on $M(1)\otimes {\mathbb C}[P]$ acting by left
multiplication with $1\otimes e\sp\lambda$. It is clear that
$$
[\omega]=e\sp{\omega_\ell}
\epsilon(\cdot,\omega_\ell),\quad\quad[\omega]\colon M(1)\otimes
{\mathbb C}[P]\to M(1)\otimes {\mathbb C}[P]
$$
is a linear bijection. It is easy to see that
\begin{equation*}
L(\Lambda_0)\overset{[\omega]}\longrightarrow
L(\Lambda_{\ell})\overset{[\omega]}\longrightarrow
L(\Lambda_{\ell-1})\overset{[\omega]}\longrightarrow \dots
\overset{[\omega]}\longrightarrow L(\Lambda_{1})
\overset{[\omega]}\longrightarrow L(\Lambda_{0}).
\end{equation*}
By using vertex operator formula (\ref{vertex operator formula})
we see that
\begin{equation}\label{the action of omega on h.w.v.}
[\omega]v_{\Lambda_0}=v_{\Lambda_\ell},\qquad
[\omega]v_{\Lambda_i}=x_{\gamma_i}(-1)v_{\Lambda_{i-1}}\quad\text{for}
\quad i=1,\dots,\ell
\end{equation}
if we properly normalize root vectors
$x_{\gamma_1},\dots,x_{\gamma_\ell}$.

By using vertex operator formula (\ref{vertex operator formula})
we see that  \ $x_\alpha
(z)[\omega]=[\omega]z\sp{\langle\omega_\ell,\alpha\rangle}x_\alpha(z)$,
or written by components
\begin{equation}\label{commutation of omega and x alpha}
x_\alpha
(n)[\omega]=[\omega]x_\alpha(n+\langle\omega_\ell,\alpha\rangle)
\quad \text{for}\quad \alpha\in R.
\end{equation}
Up to a scalar multiple, linear bijection $[\omega]$ between two
irreducible modules is uniquely determined by (\ref{commutation of
omega and x alpha}). Haisheng Li pointed out that $[\omega]\colon
L(\Lambda_i)\to L(\Lambda_{i-1})$ can be interpreted in terms of
simple currents as the identity map $\text{id\,}\colon
L(\Lambda_i)\to L(\Lambda_{i})$ if we endow the target vector space
$L(\Lambda_{i})$ the structure of $L(\Lambda_{i-1})$ with vertex
operators $Y_{L(\Lambda_{i})}(\Delta(\omega_\ell,z)\cdot,z)$, where
$\Delta(\omega,z)=z\sp\omega E\sp+(-\omega,-z)$ (cf. \cite{DLM}).

Since for fundamental modules we have $[\omega]\colon
L(\Lambda_i)\to L(\Lambda_{i-1})$, define linear bijection
$[\omega]$ on the tensor product of $k$ fundamental modules as
$$
[\omega]\otimes\dots \otimes[\omega]\colon\bigotimes_{s=1}\sp k
L(\Lambda_{i_s})\to \bigotimes_{s=1}\sp k L(\Lambda_{i_s-1}).
$$

It is clear that relation (\ref{commutation of omega and x alpha})
holds for $[\omega]=[\omega]\otimes\dots \otimes[\omega]$. In
particular,
\begin{equation}\label{commutation of omega and x}
x_\gamma (n)[\omega]=[\omega]x_\gamma(n+1) \quad \text{for}\quad
\gamma\in\Gamma.
\end{equation}
If we set $\mu\sp{+}(x_\gamma(n+1))=\mu(x_\gamma(n))$, then for
monomials relation (\ref{commutation of omega and x}) reads as
\begin{lemma}\label{commutation of omega and monomials}
\quad$x(\mu)[\omega]=[\omega]x(\mu\sp+)$.
\end{lemma}

\section{Initial conditions and simple current operator}

For an $\ell$-tuple $(a_{\ell-1},\dots,a_0)$ and the fixed level
$k$ set
$$
(a_{\ell-1},\dots,
a_1,a_0)\sp*=(a_{\ell-2},\dots,a_0,a_\ell\sp*),\quad
a_\ell\sp*=k-a_0-\dots-a_{\ell-1}.
$$
If
$$
\Lambda=k_0\Lambda_0+k_1\Lambda_1+\dots+k_{\ell-1}\Lambda_{\ell-1}+k_{\ell}\Lambda_{\ell}
$$
is dominant integral of level $k$, then $k_{\ell}\sp*=k_{\ell}$
and
$$
\Lambda\sp*=k_\ell\Lambda_0+k_0\Lambda_1+\dots+k_{\ell-2}\Lambda_{\ell-1}+k_{\ell-1}\Lambda_{\ell}
$$
is also dominant integral of level $k$. As above, we write
$v_{\mathcal K}=v_\Lambda$ and $v_{{\mathcal
K}\sp*}=v_{\Lambda\sp*}$ for ${\mathcal
K}=(k_{\ell-1},\dots,k_0)$.

Let ${\mathcal A}=(a_{\ell-1},\dots,a_0)$, where $a_i\in\mathbb
Z_+$ and $a_0+\dots+a_{\ell-1}\leq k$. Then the level $1$ initial
conditions (\ref{level 1 initial conditions}), Pauli's exclusion
principle and Dirichlet's pigeonhole principle give
\begin{align*}
&x_{\gamma_\ell}(-1)\sp{a_{\ell-1}}\dots
x_{\gamma_1}(-1)\sp{a_0}\left(v_{\Lambda_{\ell}}\sp{\otimes
a_\ell}\otimes v_{\Lambda_{\ell-1}}\sp{\otimes
a_{\ell-1}}\otimes\dots \otimes v_{\Lambda_{1}}\sp{\otimes
a_1}\otimes
v_{\Lambda_{0}}\sp{\otimes a_0}\right)\\
&=v_{\Lambda_{\ell}}\sp{\otimes
a_\ell}\otimes(x_{\gamma_\ell}(-1)v_{\Lambda_{\ell-1}})\sp{\otimes
a_{\ell-1}}\otimes\dots\otimes(x_{\gamma_2}(-1)v_{\Lambda_{1}})\sp{\otimes
a_1}\otimes(x_{\gamma_1}(-1)v_{\Lambda_{0}})\sp{\otimes a_0}\\
&=([\omega]v_{\Lambda_{0}})\sp{\otimes a_\ell\sp*}\otimes
([\omega]v_{\Lambda_{\ell}})\sp{\otimes a_{\ell-1}}\otimes\dots
\otimes([\omega]v_{\Lambda_{2}})\sp{\otimes
a_1}\otimes([\omega]v_{\Lambda_{1}})\sp{\otimes a_0}\\
&=[\omega]\left(v_{\Lambda_{0}}\sp{\otimes a_\ell\sp*}\otimes
v_{\Lambda_{\ell}}\sp{\otimes a_{\ell-1}}\otimes\dots \otimes
v_{\Lambda_{2}}\sp{\otimes a_1}\otimes v_{\Lambda_{1}}\sp{\otimes
a_0}\right),
\end{align*}
where the second equality follows from (\ref{the action of omega on
h.w.v.}). Hence we have
\begin{lemma}\label{initial conditions and omega}
\quad $x\sp{\mathcal A}v_{\mathcal A}=[\omega]v_{{\mathcal A}\sp*}$.
\end{lemma}
Lemma \ref{initial conditions and omega} and (\ref{commutation of
omega and x alpha}) imply
$$
[\omega]\colon L(\Lambda\sp*)\to L(\Lambda),\qquad [\omega]\colon
W(\Lambda\sp*)\to W(\Lambda).
$$
Note that these maps are injective. Since $L(\Lambda)$ is
irreducible $\tilde{\mathfrak g}$-module, and (\ref{commutation of
omega and x alpha}) holds, $[\omega]\colon L(\Lambda\sp*)\to
L(\Lambda)$ is a bijection.

\section{Coefficients of intertwining operators}

Let $e_1,\dots, e_\ell, e_{\ell+1}$ be the canonical basis of
${\mathbb C}\sp{\ell+1}$, viewed as a $\mathfrak g$-module, and let
$\lambda_1,\dots,\lambda_\ell,\lambda_{\ell+1}$ be the corresponding
weights of these vectors. Note that $\lambda_1=\omega_1$ and
\begin{equation}\label{x gamma is 0 on h.w.v}
x_\gamma e_j=0 \quad \text{for}\quad \gamma\in\Gamma, \
j=1,\dots,\ell.
\end{equation}
Since $\mathfrak g$-module $L(\omega_1)$ is ``on the top'' of
$\tilde{\mathfrak g}$-module $L(\Lambda_1)$, (\ref{x gamma is 0 on
h.w.v}) implies
\begin{equation}\label{x n gamma is 0 on h.w.v}
x_\gamma(n) \left(1\otimes e\sp{\lambda_j}\right)=0 \quad
\text{for}\quad \gamma\in\Gamma, \ n\geq 0, \ j=1,\dots,\ell.
\end{equation}
By (\ref{x n gamma is 0 on h.w.v}) and the commutator formula for
intertwining operators, all coefficients of $\mathcal Y(1\otimes
e\sp{\lambda_j},z)$,  $j=1,\dots,\ell$, commute with all
$x_\gamma(n)$, $\gamma\in\Gamma$, $n\in\mathbb Z$. We shall make use
of just a few of these coefficients: for $i=1,\dots,\ell$ define
$$
[i]=\text{Res}\, z\sp{-1-\langle \lambda_i,\omega_{i-1}\rangle}
c_i {\mathcal Y}(1\otimes e\sp{\lambda_i},z).
$$
By using (\ref{intertwining operators}) and (\ref{vertex operator
formula}) we see that
$$
L(\Lambda_0)\overset{[1]}\longrightarrow
L(\Lambda_{1})\overset{[2]}\longrightarrow
L(\Lambda_{2})\overset{[3]}\longrightarrow \dots
\overset{[\ell-1]}\longrightarrow L(\Lambda_{\ell-1})
\overset{[\ell]}\longrightarrow L(\Lambda_{\ell})
$$
and, with a suitable choice of $c_i$,
$$
v_{\Lambda_0}\overset{[1]}\longrightarrow
v_{\Lambda_{1}}\overset{[2]}\longrightarrow
v_{\Lambda_{2}}\overset{[3]}\longrightarrow \dots
\overset{[\ell-1]}\longrightarrow v_{\Lambda_{\ell-1}}
\overset{[\ell]}\longrightarrow v_{\Lambda_{\ell}}.
$$

Since for fundamental modules we have $[i]\colon
L(\Lambda_{i-1})\to L(\Lambda_{i})$, we may consider a linear map
$$
[i]=1\otimes\dots\otimes 1\otimes [i] \otimes 1\otimes\dots
\otimes 1
$$
on the tensor product of $k$ fundamental modules with $k_{i-1}\geq
1 $
\begin{align*}
[i]\colon & L(\Lambda_{\ell})\sp{\otimes k_\ell}\otimes\dots \otimes
L(\Lambda_{i})\sp{\otimes k_i}\otimes L(\Lambda_{i-1})\sp{\otimes
k_{i-1}}\otimes\dots \otimes L(\Lambda_{0})\sp{\otimes k_0}\\
\longrightarrow \  &L(\Lambda_{\ell})\sp{\otimes k_\ell}\otimes\dots
\otimes L(\Lambda_{i})\sp{\otimes k_i+1}\otimes
L(\Lambda_{i-1})\sp{\otimes k_{i-1}-1}\otimes\dots\otimes
L(\Lambda_{0})\sp{\otimes k_0}.
\end{align*}
Then $[i]$ maps the highest weight vector $v_\Lambda$ to a highest
weight vector
\begin{align*}
[i]\colon &  v_{\Lambda_{\ell}}\sp{\otimes k_\ell}\otimes\dots
\otimes v_{\Lambda_{i}}\sp{\otimes k_i}\otimes
v_{\Lambda_{i-1}}\sp{\otimes k_{i-1}}\otimes\dots\otimes
v_{\Lambda_{0}}\sp{\otimes k_0}\\
\mapsto \ &v_{\Lambda_{\ell}}\sp{\otimes k_\ell}\otimes\dots \otimes
v_{\Lambda_{i}}\sp{\otimes k_i+1}\otimes
v_{\Lambda_{i-1}}\sp{\otimes k_{i-1}-1}\otimes\dots\otimes
v_{\Lambda_{0}}\sp{\otimes k_0},
\end{align*}
that is
$$
[i]\colon v_{\Lambda}\mapsto v_{\Lambda'}\,,
$$
where
$$
\Lambda=k_0\Lambda_0+\dots+k_{i-1}\Lambda_{i-1}+k_i\Lambda_i+\dots+k_\ell\Lambda_\ell
$$
with $k_{i-1}\geq 1$, and
$$
\Lambda'=k_0\Lambda_0+\dots+(k_{i-1}-1)\Lambda_{i-1}+(k_i+1)\Lambda_i+\dots+k_\ell\Lambda_\ell.
$$
Note that for the corresponding $\ell$-tuples we have ${\mathcal
K}>{\mathcal K}'$. On the other hand, if ${\mathcal K}>{\mathcal
A}$, then there is a composition of operators $[i]$,
$i\in\{1,\dots,\ell\}$, such that $v_{\mathcal K}\mapsto v_{\mathcal
A}$: we apply $[1]$ $k_0-a_0$ times, then $[2]$ $k_0+k_1-a_0-a_1$
times, and so on. If we denote by $v_{\mathcal K}\longrightarrow
v_{\mathcal A}$ a composition of operators $[i]$, \
$i\in\{1,\dots,\ell\}$, then we have
\begin{lemma}\label{intertwining maps and order}
\quad $v_{\mathcal K}\longrightarrow v_{\mathcal A}$\quad if and
only if\quad ${\mathcal K}>{\mathcal A}$.
\end{lemma}

Note that on tensor products $[i]$ commutes with all $x_\gamma(n)$
for $\gamma\in\Gamma$ and $n\in\mathbb Z$ because on fundamental
modules $[i]$ commutes with all $x_\gamma(n)$ for $\gamma\in\Gamma$
and $n\in\mathbb Z$. Hence
$$
[i]\colon W(\Lambda)\to W(\Lambda').
$$

\section{Proof of linear independence}

We prove linear independence of $(k,\ell+1)$-admissible monomial
vectors
$$
x(\pi)v_{\mathcal K}= \,\dots
x_{\gamma_1}(-2)\sp{a_\ell}x_{\gamma_\ell}(-1)\sp{a_{\ell-1}}\dots
x_{\gamma_1}(-1)\sp{a_0}v_{\mathcal K}
$$
in $ W(\Lambda)$ by induction on degree
$$
n=\sum_{\gamma\in\Gamma, \, j\geq 1}j\cdot\pi(x_\gamma(-j))=1\cdot
a_0+\dots +1\cdot a_{\ell-1}+2\cdot a_\ell +\dots
$$
of monomials $x(\pi)$, considering all level $k$ modules
simultaneously.

Assume that $\sum c_\pi x(\pi)v_{\mathcal K}=0$ and that $c_\pi\neq
0$ for some $\pi$. Let ${\mathcal A}$ be a minimal $\ell$-tuple such
that
$$
c_\pi\neq 0\quad \text{for some} \quad x(\pi)=x(\mu)x\sp{\mathcal
A}=x(\mu)x_{\gamma_\ell}(-1)\sp{a_{\ell-1}}\dots
x_{\gamma_1}(-1)\sp{a_0}.
$$
Due to initial conditions we have ${\mathcal A}\leq{\mathcal K}$. If
${\mathcal A}<{\mathcal K}$, by Lemma~\ref{intertwining maps and
order} we may apply a composition of $[i]$'s for which $v_{\mathcal
K}\longrightarrow v_{\mathcal A}$ and get
$$
\sum c_\pi x(\pi)v_{\mathcal K}\longrightarrow \sum{}\sp{{}'} c_\pi
x(\pi)v_{\mathcal A}.
$$
Due to initial conditions for $v_{\mathcal A}$ and minimality of
$\mathcal A$, the sum $\sum'$ runs over all $x(\pi)$ of the form
$x(\mu)x\sp{\mathcal A}$. So by Lemma~\ref{initial conditions and
omega} and Lemma~\ref{commutation of omega and monomials} we have
$$\sum c_\pi x(\mu)x\sp{\mathcal A}v_{\mathcal A}=\sum c_\pi
x(\mu)[\omega]v_{{\mathcal A}\sp*}=[\omega]\left(\sum c_\pi
x(\mu\sp+)v_{{\mathcal A}\sp*}\right)=0,
$$
and, since the simple current operator $[\omega]$ is injective,
\begin{equation}\label{combination of smaller degree}
\sum c_\pi x(\mu\sp+)v_{{\mathcal A}\sp*}=0.
\end{equation}
The difference conditions for monomial vectors
$$
x(\mu\sp+)v_{{\mathcal A}\sp*}=\,\dots
x_{\gamma_1}(-2)\sp{a_{2\ell}}x_{\gamma_\ell}(-1)\sp{a_{2\ell-1}}\dots
x_{\gamma_1}(-1)\sp{a_\ell}v_{{\mathcal A}\sp*}
$$
are clearly satisfied, and the initial conditions for ${\mathcal
A}\sp*=(a_{\ell-2},\dots,a_0,a_\ell\sp*)$ read
\begin{align}\label{diff.cond.-init.cond.1}
&a_\ell\leq a_\ell\sp*=k-a_0-\dots-a_{\ell-1},\notag\\
&a_\ell+a_{\ell+1}\leq a_\ell\sp*+a_0=k-a_1-\dots-a_{\ell-1},\\
&\quad\vdots \notag\\
&a_\ell+\dots+a_{2\ell-1}\leq
a_\ell\sp*+a_0+\dots+a_{\ell-2}=k-a_{\ell-1}.\notag
\end{align}
We see that (\ref{diff.cond.-init.cond.1}) holds because monomial
vectors $x(\pi)v_{\mathcal K}$ satisfy difference conditions
(\ref{admissible difference conditions}), and in particular
\begin{align}\label{diff.cond.-init.cond.2}
&a_0+\dots+a_{\ell-1}+a_\ell\leq k,\notag\\
&a_1+\dots+a_{\ell-1}+a_\ell+a_{\ell+1}\leq k,\\
&\quad\vdots \notag\\
&a_{\ell-1}+a_\ell+\dots+a_{2\ell-1}\leq k.\notag
\end{align}
Hence $x(\mu\sp+)v_{{\mathcal A}\sp*}$ are $(k,\ell +1)$-admissible
vectors of degree $<n$, and our induction hypothesis together with
(\ref{combination of smaller degree}) implies that all $c_\pi=0$.
This is in contradiction with our assumption that some $c_\pi\neq
0$.

This proves that $\sum c_\pi x(\pi)v_{\mathcal K}=0$ implies $c_\pi=
0$.
%%%%%%%%%%%%%%%%%%%%%%%%%%%%%%%%%%%%%%%%%%%%%%%%%%%%%%%%%%

\end{document}